\documentclass[12pt]{amsart}

\usepackage{latexsym, amsthm, amscd, graphicx, euscript, float, verbatim, subfig}
\usepackage{amssymb}
\usepackage{amsmath}
\usepackage{graphicx}

\newcounter{figs}

\theoremstyle{remark}

\numberwithin{equation}{section}

\DeclareMathOperator{\diam}{diam}
\DeclareMathOperator{\card}{card}
\newcommand{\R}{{\mathbb R}}

\newcommand{\Rt}{{\mathbb R}^2}

\newcommand{\Rn}{{\mathbb R}^n}
\newcommand{\Rnbar}{\overline{{\mathbb R}}^n}
\newcommand{\B}{{\mathbb B}}

\newcommand{\Hn}{{\mathbb H}^n}

\newcommand{\ak}{\tilde{\alpha}}

\newcommand{\psubset}{\varsubsetneq}

\newcommand{\half}{\frac{1}{2}}

\newcommand{\PP}{ {\R^2 \setminus \{ 0 \}} }
\newcommand{\angk}[3]{\measuredangle_k (#1,#2,#3)}
\newcommand{\trik}[3]{\triangle_k (#1,#2,#3)}
\newcommand{\triks}[3]{\triangle_k^* (#1,#2,#3)}

\newcommand{\ang}[3]{\measuredangle (#1,#2,#3)}


\newcounter{minutes}
\divide\time by 60
\newcounter{hours}
\multiply\time by 60 \addtocounter{minutes}{-\time}

\begin{document}

\title[Geometry of Metrics]
{\vspace*{.1cm} Geometry of Metrics}


\author[Matti Vuorinen]{\noindent Matti Vuorinen}
\email{vuorinen@utu.fi}
\address{\newline Department of Mathematics,
\newline University of Turku,
\newline FIN-20014 Turku,
\newline Finland
}


\def\thefootnote{}
\footnotetext{
\texttt{\tiny File:~\jobname .tex,
          printed: \number\year-0\number\month-0\number\day,
          \thehours.\ifnum\theminutes<10{0}\fi\theminutes}
}


\begin{abstract}
During the past thirty years hyperbolic type metrics have become popular
tools in modern mapping theory, e.g., in the study of
quasiconformal and quasiregular maps in the euclidean $n$-space.
We study here several metrics that
one way or another are related to modern mapping theory
and point out many open problems dealing with the geometry of such metrics.
\end{abstract}

\maketitle


\section{Introduction}

Many results of
classical function theory (CFT) are more natural
when expressed in terms of the hyperbolic metric than the euclidean
metric. Naturality refers here to invariance with respect to conformal
maps or specific subgroups of M\"obius transformations. For
example, one of the corner stones of CFT, the Schwarz lemma \cite{ah},
says that an analytic function of the unit disk into itself is
a hyperbolic contraction, i.e., decreases hyperbolic distances.
Another example is Nevanlinna's principle of the hyperbolic
metric \cite[p. 50]{n}. Note that the hyperbolic metric is invariant under
conformal maps. Usual methods of
CFT such as power series, integral formulas, calculus of residues, are
mainly concerned with the local behavior of functions and do not reflect
invariance very well. The extremal length method of Ahlfors and Beurling
\cite{ab} has conformal invariance as a built-in feature and has become
a powerful tool of CFT during the past sixty years that have elapsed since
its discovery. One may even say that
conformal invariance, and thus "naturality", is one of the guiding principles
of geometric function theory.

There are serious obstacles in generalizing these ideas from the two-dimensional
case to euclidean spaces of dimension $n\ge 3\,.$ For instance, basic facts
such as multiplication of complex numbers or power series of
functions, do not make sense here. Perhaps
a more dramatic obstacle is the failure of Riemann's mapping theorem for
dimensions  $n\ge 3\,:$ according to Liouville's theorem,
conformal maps of a domain $D\subset  {\mathbb R}^n $
onto $D' \subset  {\mathbb R}^n $ are of the  form $f = g | D$ for some
M\"obius transformation $g \,.$

Here we shall study various ways to generalize the hyperbolic metric to
the $n$-dimensional case $n\ge 3\,.$ In order to circumvent the above
difficulties, we do not require complete invariance with respect
to a group of transformations, but only require
``quasi-invariance'' under transformations called quasi-isometries.
Now there are numerous
``degrees of freedom'', for instance some of the questions posed below
make sense in general metric spaces equipped with some special properties.
Therefore the problems below allow for a great number of variations, depending
on the particular metric or on the geometry of the space. The Dictionary of
Distances \cite{dd} lists hundreds of metrics.

This survey is based on my lectures held in two workshops/conferences
at IIT-Madras in December 2009 and August 2010. In December 2005 I gave
a similar survey \cite{Vu-05} and this survey partially overlaps it.
The main difference is that here mainly metric spaces are studied
while in the previous survey also  categories of maps between
metric spaces such as bilipschitz maps or quasiconformal
maps were considered. During the past decade the progress has been
rapid in this area as shown by the several recent PhD theses \cite{Ha2,he, ibr,k3,Lin, man-08, sahoo}. In fact, some of the many problems formulated in \cite{Vu-05} have been solved
 in \cite{k3,Lin, man-08}. The original informal lecture style has been
mainly kept without major changes. As in \cite{Vu-05}, several problems of varying
level of difficulty, from challenging exercises to
research problems, are given. It was assumed that the audience was
familiar with basic real and complex analysis. The interested reader
might wish to study some of the earlier surveys such as
\cite{g99, g05, v99, vu4, vu5, Vu-05}.

{\sc Acknowledgements.}  {The research of the author was supported by grant 2600066611 of the Academy of Finland.}

\section{Topological and metric spaces}

We list some basic notions from topology and metric spaces. For more
information on this topic the reader is referred to some some standard
textbook of topology such as Gamelin and Greene \cite{GG}.

The notion of a metric space was introduced by M. Frech\'et in his
thesis in 1906. It became quickly one of the key notions of topology,
functional analysis and geometry. In fact, distances and metrics occur
in practically all areas of mathematical research, see the book \cite{dd}.
Modern mapping theory in the setup of metric spaces with some additional
structure has been developed by Heinonen \cite{hei}.

\subsection{Metric space $(X,d)$}
Let $X$ be a nonempty set and let $d: X \times X \to [0,\infty)$ be
a function satisfying
\begin{enumerate}
\item[(a)] $d(x,y)= d(y,x) \,,$  for all $x,y \in X, $
\item[(b)] $d(x,y)\le d(x,z)+ d(z,y)\,,$ for all $x,y,z \in X,$
\item[(c)] $d(x,y)\ge 0 $ and $d(x,y)=0 \iff x=y\,.$
\end{enumerate}

\subsection{Examples}
\begin{enumerate}
\item $({\mathbb R}^n, |\cdot|)$ is a metric space.
\item If $(X_j,d_j), j=1,2,$ are metric spaces and $f: (X_1,d_1) \to (X_2,d_2)$
is an injection, then $m_f(x,y)= d_2(f(x), f(y))$ is a metric.
\item
If $(X,d)$ is a metric space, then also $(X, d ^{a})$ is a metric space for all
$a \in(0,1]\,.$
\item
More generally, if $h:[0,\infty)\to [0,\infty) $ is an increasing  homeomorphism
with $h(0)=0$ such that $h(t)/t$ is decreasing, then $(X,h\circ d)$ is a metric space (see e.g. \cite[7.42]{AVV}).
\end{enumerate}

\subsection{Proposition.} \label{myrho}
Let $(X,d)$ be a metric space, $0<a \le 1 \le b < \infty,$ and
$$\rho(x,y) = \max \{ d(x,y)^a, d(x,y)^b \}\,.$$
Then $$\rho(x,y) \le 2^{b-1}  (\rho(x,z)+ \rho(z,y))$$
for  all $x,y,z \in X\,.$ In particular, $\rho$ is a metric
if $b=1\,.$

\begin{proof} Fix $x,y,z \in X\,.$ Consider first the case
$d(x,y)\le 1 \,.$ Then
$$ \rho(x,y) = d(x,y)^a \le(d(x,z)+d(z,y))^a \le
d(x,z)^a+d(z,y)^a \le \rho(x,z)+ \rho(z,y) \,, $$
by an elementary inequality \cite[(1.41)]{AVV}. Next for the
case $d(x,y)\ge 1 \,$ we have
$$
\rho(x,y) = d(x,y)^b \le(d(x,z)+d(z,y))^b \le
2^{b-1}(d(x,z)^b+d(z,y)^b) \le $$
$$  2^{b-1} (\rho(x,z)+ \rho(z,y)) \,
$$
by \cite[(1.40)]{AVV}.
\end{proof}

\subsection{Remark.}
If $(X,d)$ is a metric space and $\rho$ is as defined above in \ref{myrho},
then by a result of A. H. Frink \cite{f},
there is a metric $d_1$ such that $d_1 \le \rho^{1/k}\le 4 d_1\,.$
This result was recently refined by M. Paluszynski and K. Stempak \cite{ps}.
I am indebted to J. Luukkainen for this remark.

\subsection{Uniform continuity}
Let $(X_j,d_j), j=1,2,$ be metric spaces and
$ f:(X_1,d_1) \to (X_2,d_2)$ be a continuous map. Then $f$ is uniformly
continuous (u.c.) if there exists a continuous injection $\omega:[0,t_0)\to
[0, \infty)$ such that $\omega(0)=0$ and
$$  d_2(f(x),f(y))\le \omega(d_1(x,y))\,, \mbox{ for\,\,  all\,\,} x,y\in X_1\, \mbox{ with }
d_1(x,y)<t_0.
$$

\subsection{Remarks} \label{my26}
\begin{enumerate}
\item
This definition is equivalent with the usual $(\varepsilon,\delta)$-definition \cite{GG}.
\item
If $\omega(t)=Lt$ for $t\in(0,t_0],$ then $f$ is  $L$-Lipschitz (abbr. $L$-Lip).
\item
If $\omega(t)=Lt ^a$ for some $a\in (0,1]$ and all $t\in(0,t_0],$
then $f$ is  H\"older.
\item
If $f:(X_1, d_1) \to (X_2, d_2)$ is a bijection and there is $L \ge 1$
such that $d_1(x,y)/L \le d_2(f(x), f(y)) \le L d_1(x,y)$ for all $x,y \in X_1$
then $f$ is $L$-bilipschitz. Sometimes bilipschitz maps are also called
quasi-isometries.
\item
A map is said to be an isometry if it is 1-bilipschitz.
\item
The map $f:(X,|\cdot|) \to (X,|\cdot|), X= (0,\infty), f(x)=1/x, $ for $x \in X$ is not uniformly continuous. We will later see that this map is uniformly
continuous with respect to the hyperbolic metric of $X\,.$
\item
A Lip map $h:[a,b]\to \mathbb{R}$ has a derivative a.e.
\end{enumerate}

\subsection{Balls}
Write $B_d(x_0,r)=\{ x\in X: d(x_0,x)<r  \}$ and $\overline{B}_d(x_0,r)=\{ x\in X: d(x_0,x)\le r  \}$.

\subsection{Fact} Let
$\tau = \{  B_d(x,r): x\in X, r>0 \}$ be the collection of all balls.
Then $(X,\tau \cup \{ \emptyset\} \cup \{X \})$ is a topology.

\subsection{Remarks}
\begin{enumerate}
\item We always equip a metric space with this topology.
\item The balls $\overline{B}_d(x_0,r) $ and ${B}_d(x_0,r) $ are closed
and open as point sets, resp.
\item The set $ (\mathbb Z,d), d(x,y)= |x-y|$ is a metric space. Then
$B_d(0,1)=\{0 \}$,  $\overline{B}_d(0,1)=\{-1,0,1 \}$.
Hence $clos(B_d(x_0,r))$ need not be $\overline{B}_d(x_0,r)\,.$
Also
$$ \diam(B_d(0,1))=0< \diam(\overline{B}_d(0,1))=2 \,.
$$
\item
Balls in $\mathbb{R}^n$ need not be connected (cf. below).
\item
In ${\mathbb R}^n$: balls are denoted by $B^n(x,r)$ and spheres
by $\partial B^n(x,r)=
S^{n-1}(x,r)\,.$
\end{enumerate}

\subsection{Paths}
A continuous map $\gamma:\, \Delta\to X, \Delta \subset  \mathbb R\,,$
is called a path. The length of $\gamma, \ell(\gamma)\,,$ is
$$  \ell( \gamma) = \sup \left\{ \sum_{i=1}^n d(\gamma(x_{i-1}), \gamma(x_{i}) :
\{x_0,,...,x_n\} \text{ is a subdivision of}  \, \Delta   \right \}\,.
$$
We say a path is rectifiable if $\ell(\gamma)<\infty\,.$ A rectifiable
path $\gamma:\Delta\to X$ has a parameterization in terms of
arc length $\gamma^o:[0,\ell(\gamma)] \to X $.

\subsection{Definition}
A set $G$ is connected if for all $x,y\in G$ there exists
a path $\gamma:[0,1] \to G$ such that $\gamma(0)=x,  \gamma(1)=y\,.$
Sometimes we write $\Gamma_{xy}$ for the set of all paths joining
$x$ with $y$ in $G\,.$

\subsection{Inner metric of a set $G\subset X$} For fixed
$x,y \in X$ the inner metric with respect to $G$ is defined by
$d(x,y)= \inf \{ \ell(\gamma): \gamma \in \Gamma_{xy}, \gamma \subset G  \}\,.$

\subsection{Geodesics}
A path $\gamma:[0,1]\to G$ where $G$ is a domain,
is a geodesic joining $\gamma(0)$ and $\gamma(1)$ if $\ell(\gamma)=d(\gamma(0),
 \gamma(1))$ and $d(\gamma(0),\gamma(t))+ d(\gamma(t),\gamma(1))= \ell(\gamma)$
for all $t\in (0,1).$

\subsection{Remarks}
\begin{enumerate}
\item In $(\mathbb R ^n, |\cdot|)$ the segment
$[x,y] =\{ z\in { \mathbb R} ^n : z= \lambda x+  (1-  \lambda )y,  \lambda \in [0,1] \}$ is a geodesic\,.
\item Let $G= B ^2 \setminus \{ 0 \}$ and $d$ be the inner metric of $G\,.$
There are no geodesics joining $-1/2$ and $1/2$ in $(G,d)\,.$
\end{enumerate}

\begin{figure}
\begin{center}
\includegraphics[width=2cm]{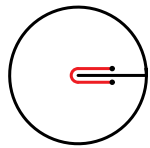}
\caption{A non-convex set which has no Euclidean geodesics \label{nonconvex}}
\end{center}
\end{figure}

\subsection{Problems} \label{convxty} (\cite[p. 322]{Vu-05})
 Let $X$ be a locally convex set in ${\mathbb R}^n$ and let $(X,d)$
be a metric space.
\begin{enumerate}
\item When are balls $B_d(x,t)$ convex for all radii $t>0$?
\item When are balls convex for small radii $t$?
\item When are the boundaries of balls nice/smooth?
\end{enumerate}

\subsection{Ball inclusion problem} \label{bip}
Suppose that $(X,d_j), j=1,2,$ determine the same euclidean topology, $X \subset
{\mathbb R}^n\,. $
Then $$B_{d_1}(x_0,r)\subset B_{d_2}(x_0,s) \subset B_{d_1}(x_0,t)$$
for some $r,s,t>0\,.$
For a fixed $s>0,$ find the best radii $r$ and $t\,.$
This problem is interesting and open for instance for
several pairs of the metrics  $d_1, d_2$ even in the special case when
 one of the metrics is the euclidean metric.

\bigskip

Note that in the above problems \ref{convxty} and  \ref{bip} we
consider the general metric space situation. It is natural to expect that
useful answers can only be given under additional hypotheses.
The reader is encouranged to find such hypotheses. In some
particular cases the problems will be studied below.

\section{Principles of geometry}

In this section we continue our list of metrics and introduce some
necessary terminology. We also outline the principles of geometry,
according to F. Klein. These principles provide a uniform view of
various geometries. In particular, the basic models of geometry:
the euclidean geometry, the geometry of the Riemann sphere and the
hyperbolic geometry of the unit ball fit into this framework.

The search for geometries leads us to compare the properties
of geometries. The Klein principles, known under the name "Erlangen Program'',
 have also paved the road for
the development of geometric function theory during the past
century. For a broad review of basic and advanced geometry
we recommend \cite{Ber} and \cite{bbi}.

\subsection{Path integrals}
For a locally rectifiable path
 $\gamma: \Delta\to X$ and a continuous function $f:\gamma\Delta \to [0,\infty] \,,$
the path integral is defined in two steps. Recall that $\gamma^o$ is the
normal representation of a rectifiable path.
\begin{enumerate}
\item[[I\!\!]]
If $\gamma$ is rectifiable, we set
$$ \int_{\gamma} f ds  = \int_0^{\ell(\gamma)}   f(\gamma^o(t)) |(\gamma^o)'(t)| \, dt \, .$$
\item[[II\!\!]]
If $\gamma$ is locally rectifiable, we set
$$ \int_{\gamma} f ds = \sup \left\{  \int_{\beta} f \, ds: \ell(\beta)<\infty,
\beta \, \mbox{ is \, a\, subpath\,of\,} \gamma \right\}\, .
$$
\end{enumerate}

\subsection{Weighted length} \label{mywlen}
Let $G\subset X$ be a domain and $w:G\to (0,\infty)$
continuous. For fixed $x,y \in D\,,$ define
$$  d_w(x,y)= \inf \{ \ell_w(\gamma) :
\gamma\in \Gamma_{xy},~\ell(\gamma)<\infty  \},
\ell_w(\gamma)= \int_{\gamma} w(\gamma(z)) \, |dz| \,.
$$
It is easy to see that $d_w$ defines a metric on $G\,$ and
$(G, d_w)$ is a metric space. If a length-minimizing
curve exists, it is called a geodesic.

\bigskip
The above construction of the weighted length \ref{mywlen} has many applications in
geometric function theory. For instance the hyperbolic and spherical
metrics are special cases of it. Our first example of \ref{mywlen} is the
quasihyperbolic metric, which has been recently studied by numerous
authors. See for instance the papers \cite{km} and \cite{KRT} in this
proceedings and their bibliographies.

\subsection{Quasihyperbolic metric}
If $w(x)=1/d(x,\partial G),$ then $d_w$ is the quasihyperbolic metric
of a domain $G \subset {\mathbb R}^n \,.$
Gehring and Osgood have proved \cite{go} that geodesics exist in this case.
Note that $w(x)=1/d(x,\partial G)$ is like a "penalty-function", the
geodesic segments try to keep away from the boundary.

\subsection{Examples}
\begin{enumerate}
\item If $G = {\mathbb H}^n = \{  x \in \mathbb R^n  :  x_n >0\}$ then
the quasihyperbolic metric coincides with the usual hyperbolic metric,
to be discussed later on.
Often the notation $\rho_{ {\mathbb H}^n}$ is used.
\item The hyperbolic metric of the unit ball  ${\mathbb B}^n$ is a weighted metric with the weight function $w(x) =2/(1- |x|^2)\,.$
Often the notation $\rho_{ {\mathbb B}^n}$ is used.


\item
In the special case when $w\equiv 1$ and the
the domain $G$ is a convex subdomain $G \subset \Rn\,,$
$d_w$ is the Euclidean distance.
The geodesics are the Euclidean segments.

\item
In the special case when $w\equiv 1$ in the non-convex set
$G=B^2 \setminus [0,1)$ geodesics do not exist (consider the points
$a=\frac{1}{2}+\frac{i}{10}$ and $b=\overline{a}$). See Fig. 1. 

\item If the construction \ref{mywlen} is applied to ${\mathbb R}^n$
with the weight function $ 1/(1+ |x|^2)$ we obtain the spherical
metric. This spaces can be identified with
the Riemann sphere
$S^{n}(e_{n+1}/2,1/2)$ equipped with the usual arc-length metric.

\item Let $X= \{ x \in {\mathbb R}:  x>0 \}$ and $w(x) = 1/x, x \in X\,.$
Then we see that $\ell_w(x,y)= |\log(x/y)|$ for all $x,y\in X\,.$
Consider again the map $f: X \to X, f(x) = 1/x, x \in X\,.$
We have seen in \ref{my26}\,(6) that it is not uniformly continuous. But it
is uniformly continuous as a map  $f: (X, \ell_w) \to (X, \ell_w) \,.$

\end{enumerate}


\subsection{The M\"obius group $\mathcal{GM}(\overline{\Rn})$}

The group of M\"obius transformations in $\overline{\Rn}$ is
generated by transformations of two types
\begin{enumerate}
\item inversions in $S^{n-1}(a,r) = \{ z\in\Rn \colon |a-z|=r \}$
$$
  x \mapsto a+\frac{r^2(x-a)}{|x-a|^2},
$$
\item reflections in hyperplane $P(a,t) = \{ x \in \Rn \colon x \cdot a = t \}$
$$
  x \mapsto x-2(x \cdot a-t)\frac{a}{|a|^2}.
$$
\end{enumerate}
If $G \subset \overline{\Rn}$ we denote by $\mathcal{GM}(G)$ the
group of all M\"obius transformations with $fG=G$.
The stereographic projection $\pi: \overline{\mathbb R}^n \to
S^n((1/2)e_{n+1}, 1/2)$ is defined by a M\"obius transformation an inversion in $S^n(e_{n+1},1)$:
$$ \pi(x) = e_{n+1}+ \frac{x-e_{n+1}}{|x-e_{n+1}|^2}\,, x \in {\mathbb R}^n, \pi(\infty)= e_{n+1}\,.$$

\subsection{Plane versus space} \label{my36}
\begin{enumerate}
\item For $n=2$ M\"obius transformations are of the form $\frac{az+b}{cz+d}, z,a,b,c,d \in {\mathbb C}\,$ with $ad-bc \neq 0\,.$

\item Recall that for $n=2$ there are many conformal maps (Riemann mapping Theorem., Schwarz-Christoffel
formula). In contrast for $n \ge 3$ conformal maps are, by
Liouville's theorem (suitable smoothness required), M\"obius
transformations.

\item Therefore conformal invariance for the space $n \ge 3$ is very different from the plane case $n=2.$
\end{enumerate}

\subsection{Chordal metric}
Stereographic projection defines the chordal distance by
$$
q(x,y) = |\pi x-\pi y| = \frac{|x-y|}{\sqrt{1+|x|^2}\sqrt{1+|y|^2}}
$$
for $x,y \in \Rnbar = {\R}^n \cup \{ \infty \}\,.$ 
\begin{figure}
\begin{center}
\vspace*{-10.5cm}
\includegraphics[width=12cm]{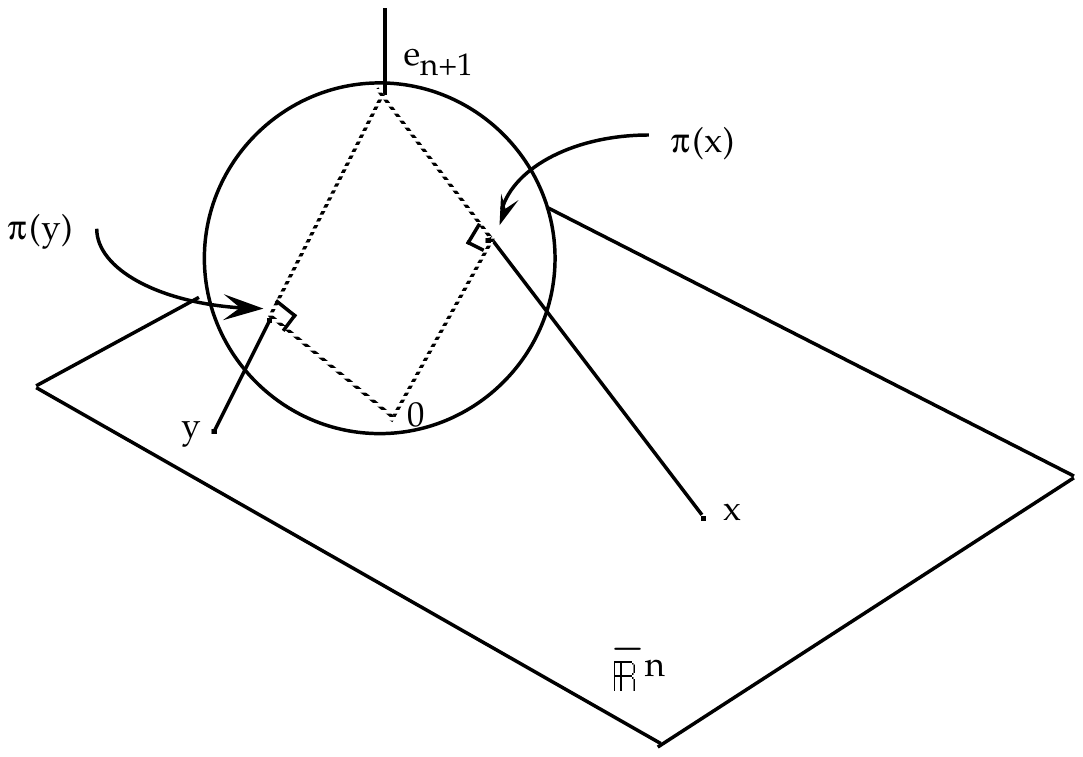}
\vspace*{-2.3cm}
 \caption{Stereographic projection \label{0701}}
\end{center}
\end{figure}
Perhaps the shortest proof of the triangle inequality for $q$
follows if we use \ref{my26}(2).

\subsection{Absolute (cross) ratio}
For distinct points $a,b,c,d \in \overline{\Rn}$ the absolute ratio
is
$$
|a,b,c,d| = \frac{q(a,c)q(b,d)}{q(a,b)q(c,d)}.
$$
The most important property is M\"obius invariance:
if $f$ is a M\"obius transformation, then $|fa,fb,fc,fd| = |a,b,c,d|$.
Permutations of $a,b,c,d$ lead to 6 different numerical values of the absolute
ratio.

\subsection{Conformal mapping}
If $G_1,G_2 \subset \Rn$ are domains and $f \colon G_1 \to G_2$ is a
diffeomorphism with $|f'(x)h| = \underbrace{|f'(x)|}_{operator\,
n.}\underbrace{|h|}_{vector\, n.}$ we call $f$ a \emph{conformal
map}. We use this also in the case $G_1,G_2 \subset \overline{\Rn}$
by excluding the two points $\{\infty, f^{-1}(\infty) \}.$

For instance, M\"obius transformations are examples of conformal maps
for all dimensions $n\ge 3$ (cf. \ref{my36} 02 )


\subsection{Linear dilatation}
Let $f : (X,d_1) \to (Y,d_2)$ be a homeomorphism, $x_0 \in X\,.$
We define the linear dilatation $H(x_0,f)$ as follows
$$  
H(x_0,f,r)   =  \frac{L_r}{l_r} \,, \quad
H(x_0,f)     =  \limsup_{r \to 0}{H(x_0,f,r)}
$$ 

\subsection{Quasiconformal maps} We adopt the definition of
V\"ais\"al\"a \cite{va1} for $K$-quasiconformal (qc) mappings. Recall that
for a $K$-qc, $K\ge 1,$  homeomorphism  $f : G\to G', G,G' \subset \R^n$
there exists a constant $H_n(K)$ such that
$\forall x_0 \in G$ $H(x_0,f) \le H_n(K) \, .$ The reader is referred to
\cite{va1} and \cite{g05} for the basic properties of quasiconformal maps.

It is well-known that conformal maps are 1-qc.
It can be also proved that $H_n(1)=1\,,$ for the somewhat tedious
details, see \cite{AVV}.


\subsection{Examples}
In most examples below, the metric spaces will have additional
structure. In particular, we will study metric spaces $(X,d) $
where the group $\Gamma$ of automorphisms of $X $ acts
transitively (i.e. given $x,y\in X$ there exists $h \in \Gamma$
such that $hx=y\,.$ We say that the metric $d$  is quasiinvariant
under the action of $\Gamma$ if there exists $C\in [1,\infty)$
such that $d(hx,hy)/d(x,y)\in [1/C,C]$ for all $x,y\in X\,, x \neq
y,$ and all $h \in \Gamma \,.$  If $C=1$, then we say that $d$ is
invariant.

\begin{enumerate}
\item
The Euclidean space $\Rn$ equipped with the usual metric
$|x-y|= (\sum_{j=1}^n (x_j-y_j)^2)^{1/2},$ $\Gamma$ is the group
of translations.

\item
The unit sphere $S^n =\{ z \in \R^{n+1} : |z|=1 \}$ equipped
with the metric of $\R^{n+1} $ and $\Gamma$ is the set of
rotations of $S^n\,.$

\end{enumerate}

\subsection{F.Klein's Erlangen Program 1872 for geometry}
\begin{itemize}
\item
use isometries ("rigid motions") to study geometry
\item
$\Gamma$ is the group of isometries
\item
two configurations are considered equivalent if they can be mapped onto each other by an element of $\Gamma$
\item
the basic "models" of geometry are %
\begin{description}
\item[(a)]
Euclidean geometry of $\R^n$
\item[(b)]
hyperbolic geometry of the unit ball $B^n$ in $\R^n$
\item[(c)]
spherical geometry (Riemann sphere)
\end{description}
\end{itemize}

The main examples of $\Gamma$ are subgroups of M\"obius transfor\-mations of
$\overline{\R}^n = \R^n \cup \{\infty\}.$


\subsection{Example: Rigid motions and invariant metrics}

        $$
        \begin{array}{c|c|l}
            X & \Gamma & \text{metric} \\
            \hline
            G & \mathcal{M}(G) & \rho_{G} \ \ \text{hyperbolic metric},
              G=B^n, \mathbb{H}^n \\
            \overline{\R}^n & \text{Iso}(\overline{\R}^n) & q \ \ \text{chordal
            metric}\\
            {\R}^n & \text{transl.} &  |\cdot| \text{Euclidean
            metric}
        \end{array}
        $$


\subsection{Beyond Erlangen, dictionary of the quasiworld}
For the purpose of studying mappings defined in subdomains of $\Rn$,
we must go beyond Erlangen, to the quasiworld, in
order to get a rich class of mappings.
\begin{center}
\begin{tabular}{rcl}
Conformal & $\to$ & ''Quasiconformal''\\
Invariance & $\to$ & ''Quasi-invariance''\\
Unit ball & $\to$ & ''Classes of domains''\\
Metric & $\to$ & ''Deformed metric''\\
World & $\to$ & ''Quasiworld''\\
Smooth & $\to$ & "Nonsmooth"\\
Hyperbolic & $\to$ & "Neohyperbolic"\\
\end{tabular}
\end{center}


\section{Classical geometries}

In this section we discuss some basic facts about the hyperbolic geometry,
already defined in Section 3 as a particular case of the weighted metric.
Some standard sources are \cite{ah2, a, be, kl}. See also \cite{Vu}.
We begin by describing the hyperbolic balls in terms of euclidean balls.
In passing we remark that this description will provide, in one concrete
case, a complete solution to the ball inclusion problem \ref{bip}.


\subsection{Comparison of metric balls}
For  $r,s>0$  we obtain the formula
\begin{eqnarray}
\rho (re_n,se_n)=\Big\vert\int_s^r \frac{dt}{t}\Big\vert
                       =\Big|\log \frac{r}{s}\Big|\;.
\label{hypform1}
\end{eqnarray}
Here $e_n=(0,...,1) \in {\mathbb R}^n\,.$
We recall the invariance property:
\begin{eqnarray}
\rho(x,y)=\rho(f(x),f(y))\;. \label{rhoinv}
\end{eqnarray}

For  $a\in \Hn$  and  $M>0$  the {\it hyperbolic ball\/}
$\{\,x\in\Hn:\rho(a,x) <M\,\}$  is denoted by  $D(a,M)$. It is
well known that  $D(a,M)=B^n(z,r)$ for some  $z$ and $r$ (this
also follows from (\ref{rhoinv})!\thinspace). This fact together
with the observation that  $\lambda te_n,\,(t/\lambda)e_n \in
\partial D(te_n,M)$, $\lambda=e^M$ (cf.\ (\ref{hypform1})), yields
\begin{equation}
\left\{\begin{array}{ll}
D(te_n,M)=B^n\bigl((t\cosh M)e_n,t\sinh M\bigr)\;,&{}\\
B^n(te_n,rt)\subset D(te_n,M)\subset B^n(te_n,Rt)\;,&{}\\
r=1-e^{-M}\;,\;\;R=e^M-1\;.&{}
\end{array}\right.
\label{hypballs}
\end{equation}

\begin{figure}
\begin{center}
\includegraphics[width=9cm]{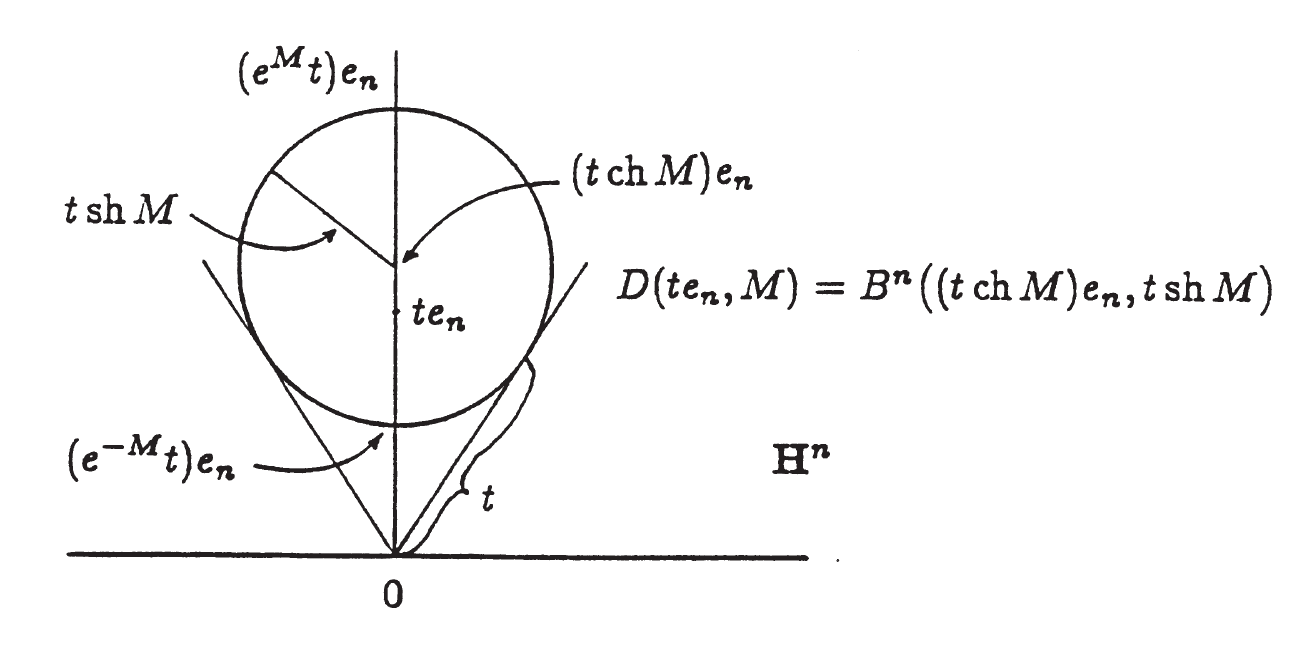}
\end{center}
\caption{The hyperbolic ball $D(t e_n,M)$ as a Euclidean ball.}
\end{figure}

\medskip
It is well known that the balls  $D(z,M)$  of  $(\B^n,\rho)$  are
balls in the Euclidean geometry as well, i.e. $D(z,M)=B^n(y,r)$
for some  $y\in \B^n$ and  $r>0$. Making use of this fact, we shall
find  $y$ and $r$. Let  $L_z$  be a Euclidean line through  $0$
and $z$ and  $\{z_1,z_2\}=L_z\cap \partial D(z,M)$, $|z_1|\le
|z_2|$. We may assume that  $z\ne 0$ since with obvious changes
the following argument works for  $z=0$ as well. Let  $e=z/|z|$
and  $z_1=se$, $z_2=ue$, $u\in (0,1)$, $s\in (-u,u)$. Then it
follows that
\begin{eqnarray*}
\rho(z_1,z)&=&\log\Bigl(\frac{1+|z|}{1-|z|}\cdot\frac{1-s}{1+s}\Bigr)=M\;,\\
\rho(z_2,z)&=&\log\Bigl(\frac{1+u}{1-u}\cdot\frac{1-|z|}{1+|z|}\Bigr)=M\;
\end{eqnarray*}

Solving these for $s$ and $u$ and using the fact that
$$D(z,M)=B^n\bigl({\textstyle \frac{1}{2}} (z_1+z_2),\,
{\textstyle \frac{1}{2}}|u-s|\bigr)$$ one obtains the following
formulae:
\begin{equation}\label{2.22}
\left\{\begin{array}{l}
D(x,M)=B^n(y,r)\\
{}\\
{\displaystyle y=\frac{x(1-t^2)}{1-|x|^2t^2}\;,\;\;
r=\frac{(1-|x|^2)t}{1-|x|^2t^2}}\;,\;\;t=\tanh {\textstyle
\frac{1}{2}} M\;,
\end{array}\right.
\end{equation}
and
\begin{equation*}
\left\{\begin{array}{l} B^n\bigl(x,\,a(1-|x|)\bigr)\subset D(x,M)
    \subset B^n\bigl(x,\,A(1-|x|)\bigr)\;\,,\\
{}\\
{\displaystyle  a=\frac{t(1+|x|)}{1+|x| t}\;,\;\;
A=\frac{t(1+|x|)}{1-|x| t}}\;,\;\;t=\tanh {\textstyle \frac{1}{2}}
M\;.
\end{array}\right.
\end{equation*}

A special case of (\ref{2.22}):
\begin{equation*}
D(0,M)=B_{\rho}(0,M)=B^n(\tanh \half M)\;.
\end{equation*}

For a given pair of points $x,y \in {\mathbb R}^n$ and a number $t>0\,,$
an Apollonian sphere is the set of all points $z$ such that $|z-x|/|z-y|=t \,.$
It is easy to show that, given $x \in B^n\,,$ hyperbolic spheres with hyperbolic
center $x$ are Apollonian spheres w.r.t. the points $x, x/|x|^2,$ see \cite{kv}.

Note that balls in chordal metric can be similarly described in terms of the
euclidean balls, see \cite{AVV}.


{
\subsection{Hyperbolic metric of the unit ball $B^n$} \label{my42}
{\bf Four equivalent definitions of the hyperbolic metric $\rho_{B^n}\,.$ }
\begin{enumerate}
\item $\rho_{B^n} = m_w$, $w(x) = \frac{2}{1-|x|^2}\,.$
\item $\textnormal{sinh}^2 \frac{\rho_{B^n}(x,y)}{2} = \frac{|x-y|^2}{(1-|x|^2)(1-|y|^2)}\,.$

\item $\rho_{B^n}(x,y) = \sup \{ \log |a,x,y,d| \colon a,d \in \partial B^n \}\,.$
\item $\rho_{B^n}(x,y) = \log |x_*,x,y,y_*|\,.$
\end{enumerate}
}
The hyperbolic metric is invariant under the action of $\mathcal{GM}(B^n)$,
i.e. $\rho(x,y)=\rho(h(x),h(y))$ for all $x,y \in B^n$ and all
$h \in \mathcal{GM}(B^n)\,.
$

\subsection{The hyperbolic line through $x,y$}
The hyperbolic geodesics between $x,y$ in the unit ball are the circular arcs
joining $x$ and $y$ orthogonal to $\partial B^n$. 

  \begin{center}
   \begin{figure}
\includegraphics[width=0.2\textwidth]{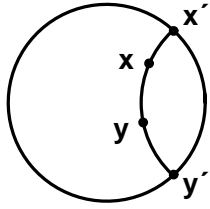}
     \caption{Hyperbolic lines are circular arcs perpendicular
     to $\partial B^n$ and $\rho_{B^n}(x,y)= \log |x',x,y,y'|.$
\label{hypebolicgeodesic2}}
   \end{figure}
\vspace{.51cm}
  \end{center}


\vspace{ -1cm}
\begin{figure}
\begin{center}
  \includegraphics[width=0.75\textwidth]{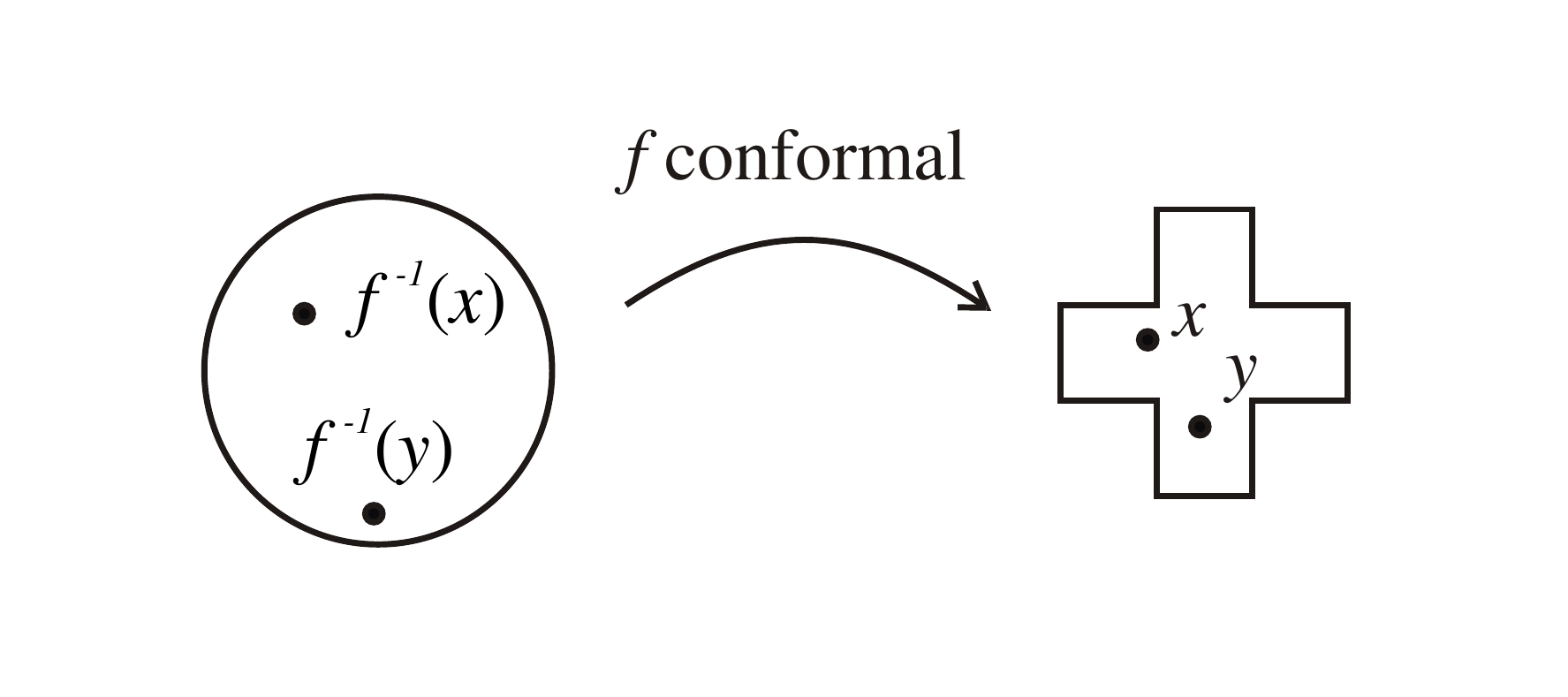}
  \end{center}
\vspace{ -1cm}
\caption{Definition of $\rho(x,y)$ in terms of
$\rho_{B^n}(f^{-1}(x),f^{-1}(y)),$ $f$ conformal.}
\end{figure}

\bigskip
\bigskip

\subsection{Hyperbolic metric of $G=f(B^2)$, $f$ conformal}
In the case where $G_k = f_k(B^2)$ and $f_k$ is conformal, $k=1,2\,,$
it follows that
if $h \colon G_1 \to G_2$ is conformal, then the hyperbolic metric
is invariant under $h\,,$ i.e., $\rho_{G_1}(x,y) =
\rho_{G_2}(hx,hy)$. Thus we may use explicit conformal maps to evaluate
the hyperbolic metrics in cases where such a map is known.


\begin{figure}
\vspace{ -1cm}
\begin{center}
  \includegraphics[width=0.75\textwidth]{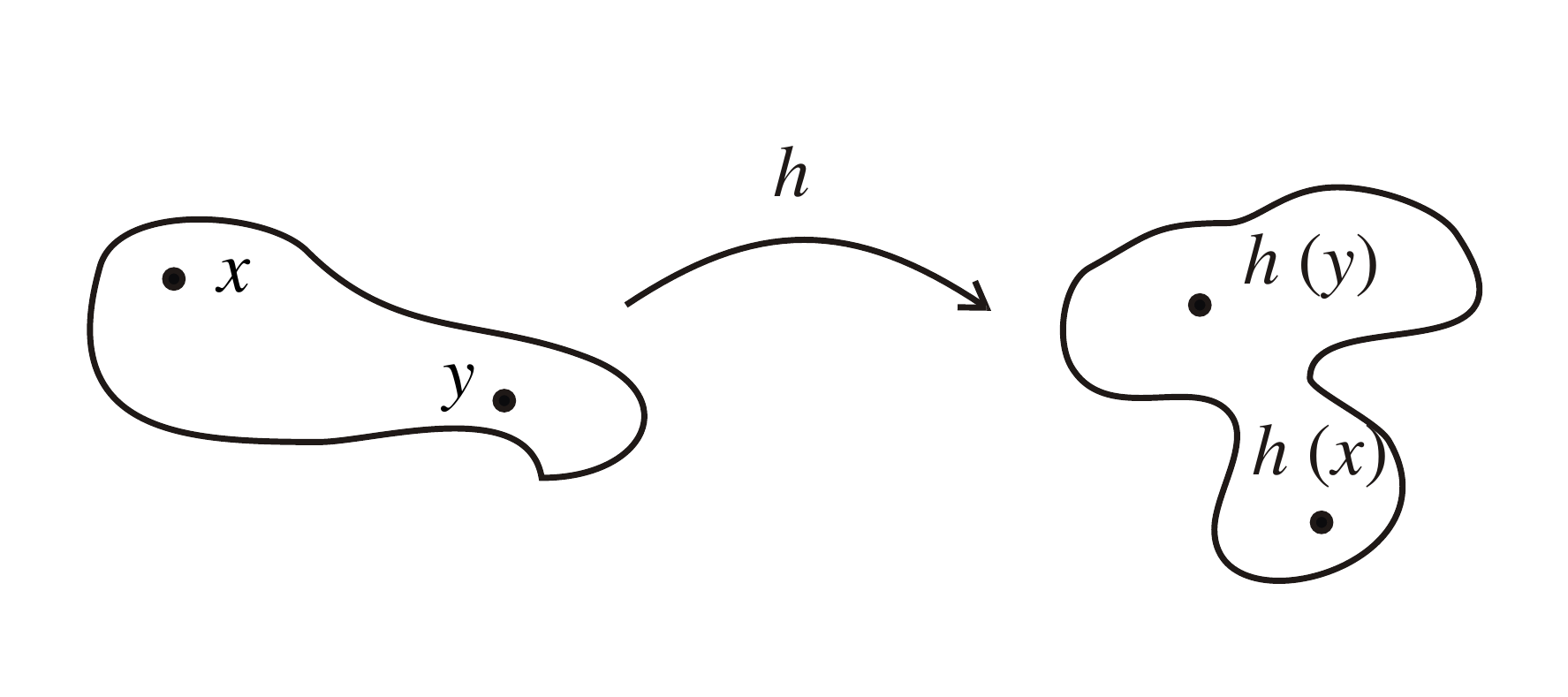}
\end{center}
\caption{Invariance of the hyperbolic metric under conformal map.}
\vspace{ -.5cm}
\end{figure}
\bigskip
For $n=2$ one can generalize the hyperbolic metric, using covering
transformations, to a domain $G \subset \overline{\R}^2$ with $\card
(\overline{\R}^2 \setminus G) \ge 3$ \cite{kl}.

The formula for the hyperbolic metric of the unit ball given by \ref{my42}(2)
is relatively complicated. Therefore various comparison functions have been
introduced. We will now discuss two of them.

\subsection{The distance ratio metric $j_G$}
For $x,y\in G$ the {\em distance ratio metric $j_G$} is defined \cite{Vu2} by
$$ j_G(x,y)=\log\left(1+\frac{|x-y|}{\min\{d(x),d(y)\}}\right)\,.
$$
The inequality $j_G\le \delta_G\le \tilde{j}_G\le 2j_G$ holds for every open set $G\psubset \Rn$,
where the metric $\tilde{j}_G$ (cf. \cite{go}) is a metric defined by
$$\tilde{j}_G(x,y)=\log\left(1+\frac{|x-y|}{d(x)}\right)\left(1+\frac{|x-y|}{d(y)}\right)\,.
$$

We collect the following well-known facts:
\begin{enumerate}
 \item Inner metric of the $j_G$ metric is the quasihyperbolic metric $k_G \,.$
\item $k_G(x,y) \le 2j_G(x,y)$ for all $x,y\in G$ with $ j_G(x,y)< \log(3/2)\,,$
see \cite[3.7(2)]{Vu}.
\item Both $k_G$ and $j_G$ define the Euclidean topology.
\item $j_G$ is not geodesic; the balls $B_j(z,M)= \{x \in G: j_G(z,x) <M\}$
may be disconnected for large $M\,.$
\end{enumerate}

If we compare the density functions of the hyperbolic and the quasihyperbolic metrics of
$B^n$, it will lead to the observations that
\begin{equation} \label{k-and-rho}
\rho_{B^n}(x,y)/2 \le k_{B^n}(x,y)\le \rho_{B^n}(x,y)
\end{equation}
for all $x,y\in B^n\,.$

The following proposition gathers together several
basic properties of the metrics $k_G$ and $j_G$, see for instance \cite{GP,Vu}.

\subsection{Proposition}\label{k-equals-j} (\cite{KSV})
\begin{enumerate}
\item \label{j-le-k} For a domain $G\subset {\mathbb R}^n, x,y \in G,$
and with $L= \inf \{ \ell(\gamma) : \gamma\in \Gamma_{x,y} \}\,,$ we have
$$k_G(x,y) \ge \log\left(1+ \frac{L}{\min\{\delta(x), \delta(y)\}}\right) \ge j_G(x,y)\,.$$
\item\label{ksv-sub1} For $x\in B^n$ we have
$$k_{B^n}(0,x)=j_{B^n}(0,x)=\log \frac{1}{1-|x|}\,.
$$
\item\label{ksv-sub2} Moreover, for $b\in S^{n-1}$ and $0<r<s<1$ we have
$$k_{B^n}(br,bs)=j_{B^n}(br,bs)=\log \frac{1-r}{1-s}\,.
$$
\item\label{ksv-sub3} Let $G\psubset\Rn$ be any domain and $z_0\in G$.
Let $z\in\partial G$ be such that $\delta(z_0)=|z_0-z|$. Then for all
$u,v\in [z_0,z]$ we have
$$k_G(u,v)=j_G(u,v)= \left|\log
\frac{\delta(z_0)-|z_0-u|}{\delta(z_0)-|z_0-v|}\right|
=\left|\log \frac{\delta(u)}{\delta(v)}\right|.$$
\item\label{rhoj} For $x,y\in B^n$ we have
$$ j_{B^n}(x,y)\le \rho_{B^n}(x,y) \le 2 j_{B^n}(x,y)
$$
with equality on the right hand side when $x=-y\,.$
\item\label{jk}
For $0<s<1$ and $x,y\in B^n(s)$ we have
$$ j_{B^n}(x,y) \le k_{B^n}(x,y) \le (1+s)\, j_{B^n}(x,y).
$$
\end{enumerate}
\begin{proof} (1) Without loss of generality we may assume that
$\delta(x)\le \delta(y).$ Fix $\gamma \in \Gamma(x,y)$ with
arc length parameterization $\gamma:[0,u]\to G, \gamma(0) =x, \gamma(u)=y$
\begin{eqnarray*} \ell_k(\gamma) = \int_{0}^u\frac{|\gamma'(t)|\,dt}{d(\gamma(t), \partial G)}\ge
\int_0^u \frac{dt}{\delta(x)+t}
 & = & \log \frac{\delta(x)+u}{\delta(x)}\\
& \ge & \log\left(1+ \frac{|x-y|}{\delta(x)}\right)=
j_G(x,y)\,.
\end{eqnarray*}

(2) We see by (\ref{j-le-k}) that
$$ j_{B^n}(0,x) = \log\frac{1}{1-|x|} \le k_{B^n}(0,x)
\le \int_{[0,x]} \frac{|dz|}{\delta(z)}= \log\frac{1}{1-|x|}
$$
and hence $[0,x]$ is the $k_{B^n}$-geodesic
between $0$ and $x$ and the equality in (\ref{ksv-sub1}) holds.

The proof of (\ref{ksv-sub2}) follows from (\ref{ksv-sub1}) because
the quasihyperbolic length is additive along a geodesic
$$ k_{B^n}(0,bs) = k_{B^n}(0,br)+ k_{B^n}(br,bs)\,.
$$

The proof of (\ref{ksv-sub3}) follows from (\ref{ksv-sub2}).

The proof of (\ref{rhoj}) is given in \cite[Lemma 7.56]{AVV}.

For the proof of the last statement see \cite{KSV}.
\end{proof}

In view of (\ref{k-and-rho}) and Proposition \ref{k-equals-j} we see
that for the case of the unit ball, the metrics $j,k, \rho$ are
closely related.


\section{Metrics in particular domains:
uniform, quasidisks}

We have seen above that for the case of the unit ball, several metrics
are equivalent. This leads to the general question: Suppose that given
a domain $G \subset {\mathbb R}^n$ we have two metrics $d_1, d_2$ on
$G\,.$ Can we characterize those domains $G\,,$
which for a fixed constant $c > 0$ satisfy
 $d_1(x,y) \le c d_2(x,y)$ for all $x,y \in G\,.$ As far as
we know, this is a largely open problem. However, the class of domains
characterized by the property that the quasihyperbolic and the
distance ratio metric have a bounded quotient, coincides with
the very widely known
class of uniform domains introduced by Martio and Sarvas \cite{MS}.

There is more general class of domains, so called $\varphi$-uniform
domains, which contain the uniform domains as special case which
we will briefly discuss \cite{Vu2}.

It is easy to see that for a general domain the quasihyperbolic and
distance ratio metrics both define the euclidean topology, in fact
we can solve the ball inclusion problem \ref{bip} easily, see \cite[(3.9)]{Vu}
for the case of the quasihyperbolic metric. Although some progress has been
made on this problem recently in \cite{kv2}, the problem is not completely
solved in the case of metrics considered in this survey.

\subsection{Uniform domains and constant of uniformity}
The following form of the definition of the uniform domain is due to
Gehring and Osgood \cite{go}. As a matter of fact, in \cite{go} there
was an additive constant in the inequality (\ref{unif-GO-eq}), but
it was shown in \cite[2.50(2)]{Vu2} that the constant
can be chosen to be $0\,.$

\subsection{Definition}\label{unif-GO}
A domain $G\psubset \Rn$ is called uniform, if there exists a number $A\ge 1$ such that
\begin{equation}\label{unif-GO-eq}
k_G(x,y)\le A\,j_G(x,y)
\end{equation}
for all $x,y\in G$. Furthermore, the best possible number
$$A_G:=\inf\{A\ge 1:\, A \mbox{ satisfies } (\ref{unif-GO-eq})\}
$$
is called the uniformity constant of $G$.

Our next goal is to explore domains for which the uniformity constant
can be evaluated or at least estimated. For that purpose we consider
some simple domains.

\subsection{ Examples of quasihyperbolic geodesics } \label{my51}
\begin{enumerate}
\item
For the domain $\Rn\setminus\{0\}$ Martin and Osgood (see \cite{MO}) have
determined the geodesics. Their result states that given $x,y\in \Rn\setminus\{0\}$,
the geodesic segment can be obtained as follows: let $\varphi$ be the angle
between the segments $[0,x]$ and $[0,y]$, $0<\varphi<\pi$. The triple
$0,x,y$ clearly determines a $2$-dimensional plane $\Sigma$, and the geodesic
segment connecting $x$ to $y$ is the logarithmic spiral in $\Sigma$ with polar
equation
$$r(\omega)=|x|\,\exp \left(\frac{\omega}{\varphi}\log \frac{|y|}{|x|}\right)\,.
$$
In the punctured space the quasihyperbolic distance is given by the formula
$$k_{\Rn\setminus\{0\}}(x,y)=\sqrt{\varphi^2+\log^2\frac{|x|}{|y|}}\,.
$$
\item \cite{Lin}
Let $\varphi\in (0,\pi]$ and $x,y\in S_\varphi=\{(r,\theta)\in \Rt:0<\theta<\varphi\}$, the angular domain.
Then the quasihyperbolic geodesic segment is a curve consisting of line segments and
circular arcs orthogonal to the boundary. If $\varphi\in (\pi,2\pi)$, then the
geodesic segment is a curve consisting of pieces of three types:
line segments, arcs of logarithmic spirals and
circular arcs orthogonal to the boundary.
\item \cite{Lin}
In the punctured ball $\B^n\setminus\{0\}$, the quasihyperbolic geodesic segment is a curve
consisting of arcs of logarithmic spirals and geodesic segments of the quasihyperbolic metric
of $B^n$.

\end{enumerate}

\begin{figure} \protect{\label{mykperj}}
\vspace*{-4.5cm}
\includegraphics[width=11cm]{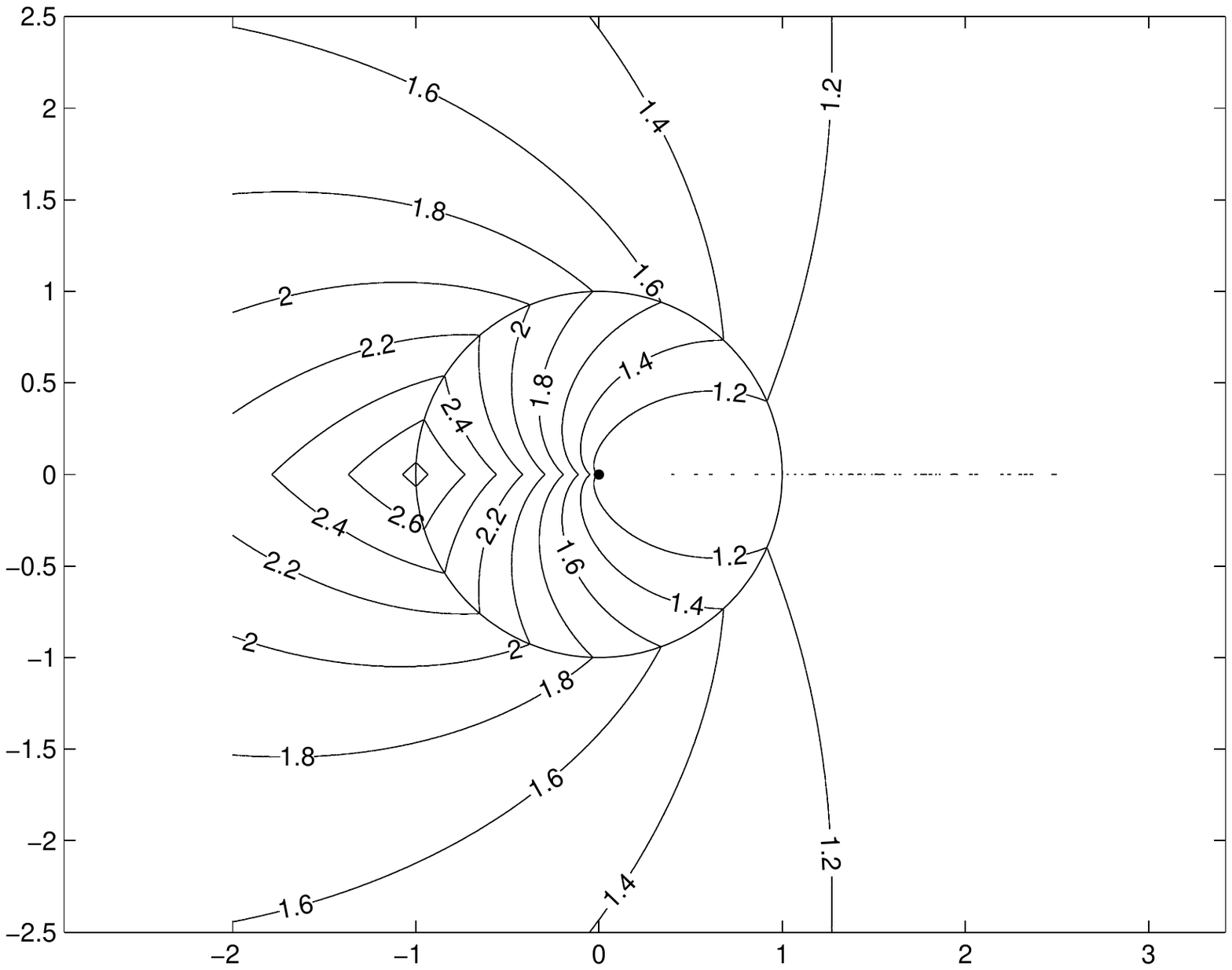}
\vspace*{-3.5cm}
\caption{Sets $\{  z :  k_G(1,z)/j_G(1,z) = c\}\,.$ }
\end{figure}


The above formula \ref{my51}(1) shows that the quasihyperbolic metric
of $G= \Rn\setminus\{0\}$  is
invariant under the inversion $x \mapsto x/|x|^2 \,$ which maps $G$ onto
itself. It is also easy to see that for this domain $G$ also  $j_B$ has the
same invariance property. Next, for this domain $G$ and a given number $c>1\,,$
the sets $\{  x :  k_G(1,x)/j_G(1,x) = c\}\,$ are illustrated. The invariance
under the inversion is quite apparent. The same formula 5.3 (1) is also
discussed in \cite{km}.

\bigskip

We now give a list of constants of uniformity for a few specific domains
following H. Lind\'en \cite{Lin}.
\begin{enumerate}
\item
For the domain $\Rn\setminus\{0\}$, the uniformity constant is given by
(cf. Figure 7)
$$A_{\Rn\setminus\{0\}}=\pi/\log 3\approx 2.8596\,.
$$
\item
Constant of uniformity in the punctured ball $B^n\setminus\{0\}$
is same as that in $\Rn\setminus\{0\}$.
\item
For the angular domain $S_\varphi$, the uniformity constant is given by
$$A_{S_\varphi}=\frac{1}{\sin \frac{\varphi}{2}}+1
$$
when $\varphi\in (0,\pi]$.
\end{enumerate}

There are numerous characterizations of quasidisks, i.e. quasiconformal
images of the unit disk under a quasiconformal map. E.g. it is known that
a simply connected domain is a quasidisk if and only if it is a uniform
domain, see \cite{g99}.

\subsection{$\varphi$-uniform domains (\cite{Vu2}).} Let
$\varphi:[0,\infty) \to [0,\infty)$ be a homeomorphism. We say that a
domain $G \subset {\mathbb R}^n$ is $\varphi$-uniform if
$$
k_G(x,y) \le \varphi( |x-y|/\min \{d(x, \partial G), d(y, \partial G) \})
$$
holds for all $x,y \in G\,.$

\bigskip

In \cite{Vu2} $\varphi$-uniform domains were introduced for the purpose
of finding a wide class of domains where various conformal invariants
could be compared to each other. Obviously, uniform domains form a subclass.
Recently, many examples of these domains were given in \cite{KSV}.
This class of domains is relative little investigated and there are many
interesting questions even in the case of plane simply connected
$\varphi$-uniform domains. This class of plane domains contains e.g.
all quasicircles. Because for a quasicircle $C$ the both components of
${\mathbb C} \setminus C$ are quasidisks, we could ask the following
question. Suppose that $C$ is a Jordan curve in the plane dividing
thus ${\mathbb C} \setminus C$ into two components, one of which is
a $\varphi$-uniform domain. Is it true that also the other component is
a $\varphi_1$-uniform domain for some function $\varphi_1$? This
question was recently answered in the negative in \cite{HKSV}.

\subsection{Open problem.} Is it true that there are
 $\varphi$-uniform domains $G$ in the plane such that the Hausdorff-dimension
of $\partial G$ is two?

\bigskip
Recall that for quasicircles this is not possible by \cite{gv}. P. Koskela
has informed the author that Tomi Nieminen has done some work on this
problem.

\section{Hyperbolic type geometries}

In this section we discuss briefly two metrics, the Apollonian metric
$\alpha_G$ and a M\"obius invariant metric $\delta_G$ introduced by
P. Seittenranta \cite{Se} and formulate a few open problems. For
the case of the unit ball, both metrics coincide with the hyperbolic
metric. For other domains they are quite different: while
$\delta_G$ is always a metric, for domains with small boundary $\alpha_G$
may only be a pseudometric. The Apollonian metric was introduced in 1934
by D. Barbilian \cite{barb, bs}, but forgotten for many years.
A. Beardon \cite{Be} rediscovered it independently in 1998
and thereafter it has been studied very intensively by many authors:
see, e.g., Z. Ibragimov \cite{ibr},
P. H\"ast\"o \cite{Ha2, Ha5, Ha-05, Ha5, hi, HPS,HKSV},
S. Ponnusamy \cite{HPWS, HPWW}, S. Sahoo \cite{sahoo}.
See also D. Herron, W. Ma and D. Minda \cite{hmm}.

\subsection{Apollonian metric of $G \psubset \Rn$}
$$
\alpha_G(x,y) = \sup \{ \log |a,x,y,b| \colon a,b \in \partial G \}.
$$
\begin{itemize}
\item $\alpha_G$ agrees with $\rho_G$, if $G$ equals $B^n$ and $H^n$.
\item $\alpha_{hG}(hx,hy) = \alpha_G(x,y)$ for $h \in \mathcal{GM}(\Rn)$
\item $\alpha_G$ is a pseudometric if $\partial G$ is ''degenerate''
\end{itemize}

\subsubsection{Facts}
\begin{enumerate}
\item The well-known sharp relations $\alpha_G\le 2j_G$ and $\alpha_G\le 2 k_G$ are due to
Beardon \cite{Be}.
\item $\alpha_G$ does not have geodesics.
\item The inner metric of the Apollonian metric is called the Apollonian inner metric and
it is denoted by $\ak_G$ (see \cite{Ha2,Ha5,HPS}).
\item We have $\alpha_G\le \ak_G\le 2k_G$.
\item $\ak_G$-geodesic exists between any pair of points in $G\psubset \Rn$ if
$G^c$ is not contained in a hyperplane \cite{Ha5}.
\end{enumerate}

\begin{figure}
\begin{center}
\vspace{-.5cm}
\includegraphics[width=7cm]{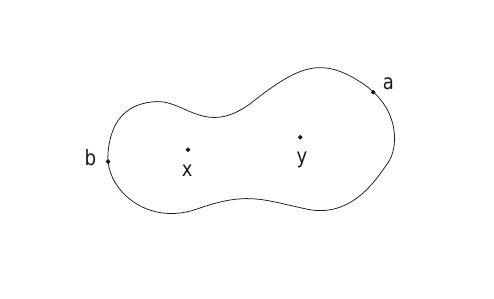} 
\caption{A quadruple of points admissible for the definition of
the Apollonian metric.}
\vspace{-.5cm}  
\end{center}
\end{figure}

\subsection{A M\"obius invariant metric $\delta_G$} For $x,y\in G\psubset \Rn$,
Seittenranta \cite{Se} introduced the following metric
$$\delta_G(x,y)=\sup_{a,b\in\partial G}\log \{1+|a,x,b,y|\}\,.
$$

\subsubsection{Facts} \cite{Se}
\begin{enumerate}
\item
The function $\delta_G$ is a metric.
\item
$\delta_G$ agrees with $\rho_G$, if $G$ equals $B^n$ or $H^n$.
\item It follows from the definitions that $\delta_{\Rn\setminus\{a\}}=j_{\Rn\setminus\{a\}}$ for all
$a\in \Rn$.
\item
$\alpha_G\le \delta_G\le \log (e^{\alpha_G}+2)\le \alpha_G+3$. The first two inequalities are
best possible for $\delta_G$ in terms of $\alpha_G$ only \cite{Se}.
\end{enumerate}


\subsection{Open problem}
Define
$$ m_{B^n}(x,y):= 2\log\left(1+\frac{|x-y|}{2\min\{d(x),d(y)\}}\right)\, .
$$

Then $m_{B^n}(x,y)$ is not a metric. In fact, any choice of the points
on a radial segment will violate the triangle inequality. It is easy to see that
 $k_{B^n}(x,-x) =m_{B^n}(x, -x)$. We do not know whether
$k_{B^n}(x,y) \le m_{B^n}(x, y)$ for all $x,y \in B^n\,.$
 If the inequality  holds, then certainly
$k_{B^n}\le 2m_{B^n}\le 2j_{B^n}$.


\subsection{Diameter problems}
There exists a domain $G\psubset\Rn$ and $x\in G$ such that $j(\partial B_j(x,M))\neq 2M$
for all $M>0$. Indeed, let $G=B^n$. Choose $x\in (0,e_1)$ and consider the $j$-sphere
$\partial B_j(0,M)$ for $M=j_G(x,0)$. Now, $B_j(0,M)$ is a Euclidean ball with radius
$|x|=1-e^{-M}$. The diameter of the $j$-sphere $\partial B_j(0,M)$ is
$$j_G(x,-x)=\log\left(1+\frac{|2x|}{d(x)}\right)=\log\left(1+\frac{2-2e^{-M}}{e^{-M}}\right)
=\log (2e^M-1)\,.
$$
We are interested in knowing whether $j_G(x,-x)=2M$ holds,
equivalently in this case,
$(e^M-1)^2=0$ which is not true for any $M>0$. Therefore, we always have $j_G(x,-x)<2M$
and the diameter of $\partial B_j(0,M)$ is less than twice the radius $M$.
Note that there is no geodesic of the $j_G$ metric joining $x$ and $-x\,.$

For a convex domain $G$, it is known by Martio and V\"ais\"al\"a \cite{MV} that
$k(\partial B_k(x,M))=2M$. However, we have the following open problem.

\subsection{Open problem}
Does there exist a number $M_0>0$ such that for all $M\in (0,M_0]$ we have $k(\partial B_k(x,M))=2M$. For the case of plane domains, this problem was
studied by Beardon and Minda \cite{bm2}.

\subsection{Convexity problem \cite{Vu-05}}

Fix a domain $G \subsetneq \Rn$ and neohyperbolic metric
$m $ in a collection of metrics (e.g. quasihyperbolic, Apollonian, $j_G$, hyperbolic metric of
a plane domain etc.). Does there exist constant $T_0 > 0$ such that the ball
$B_m(x,T) = \{ z \in G \colon m(x,z) < T \}$,
is convex (in Euclidean geometry) for all $T \in (0,T_0)$?


\subsection{Theorem}[\cite{k1}]
For a domain $G \subsetneq \Rn$ and $x \in G$ the $j$-balls
$B_j(x,M)$ are convex if and only if $M \in (0,\log 2]$.

\begin{figure}
\includegraphics[width=3cm]{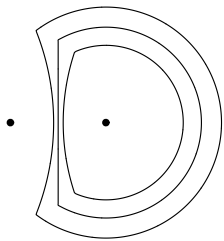}
\caption{Boundaries (nonsmooth!) of $j$-disks $B_{j_{\R^2 \setminus \{ 0 \}}}(x,M)$
with radii $M=-0.1+\log 2$, $M=\log 2$ and $M=0.1+\log 2$.}
\end{figure}

%

\subsection{Theorem}[\cite{k2}, \cite{MO}]
For $x \in {\mathbb R}^2 \setminus \{ 0 \}$
the quasihyperbolic disk $B_k(x,M)$ is
strictly convex iff $M \in (0,1]\,.$


\begin{figure}
\includegraphics[width=3cm]{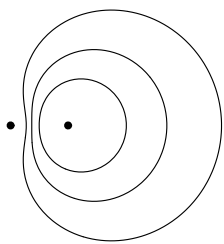}
\caption{Boundaries of quasihyperbolic disks $B_{k_{\R^2 \setminus \{ 0 \}}}(x,M)$
with radii $M=0.7$, $M=1.0$ and $M=1.3$.}
\end{figure}

Some of the convexity results of Kl\'en have been extended to Banach
spaces by A. Rasila and J. Talponen \cite{rt}. See also \cite{KRT}.


\newcommand{\Ball}[3]{B_{#1}(#2,#3)}
\newcommand{\N}{\mathbb{N}}
\newcommand{\e}{\mathbf{e}}

If a metric space is geodesic, then all metric balls are connected. For nongeodesic metric spaces the connectivity of metric balls depends on the setting.
For example, chordal balls are always connected while
$j$-balls need not be connected \cite[Remark 4.9 (2)]{k1}.
See also \cite{KRT}.

\subsection{Lemma}\label{sizeofcomponent}
  Let $G \subset \Rn$ be a domain, $x \in G$, and $r > 0$. Then for each connected component $D$ of $\Ball{j}{x}{r}$ we have
  \[
    \diam_k(  D ) \le c(r,n).
  \]

\section{Complement of the origin}

We have already seen that the quasihyperbolic metric has a simple
formula for the complement of the origin. Even more is true: many
results of elementary plane geometry hold, possibly with minor
modifications, in the quasihyperbolic geometry. Geometrically
we can view $(G, k_G ), G= {\mathbb R}^2 \setminus \{ 0 \}\,$ as
a cylindrical surface embedded in ${\mathbb R}^3$, cf. \cite{k3}.

Therefore many basic results of euclidean geometry hold for
$(G, k_G )$ as such or with minor modifications. Some of these
results are listed below.



\subsection{Theorem}[Law of Cosines]\label{mainthm} (\cite{k3})
  Let $x,y,z \in \PP$.

  \noindent(i) For the quasihyperbolic triangle $\trik{x}{y}{z}$
  $$
    k(x,y)^2 = k(x,z)^2+k(y,z)^2-2k(x,z)k(y,z)\cos \angk{y}{z}{x}.
  $$

  \noindent(ii) For the quasihyperbolic trigon $\triks{x}{y}{z}$
  $$
    k(x,y)^2 = k(x,z)^2+k(y,z)^2-2 k(y,z)k(z,x) \cos \angk{y}{z}{x}-4\pi(\pi-\alpha),
  $$
  where $\alpha = \ang{x}{0}{y}$.


\begin{figure}[htp]
  \begin{center}
    \includegraphics[height=5cm]{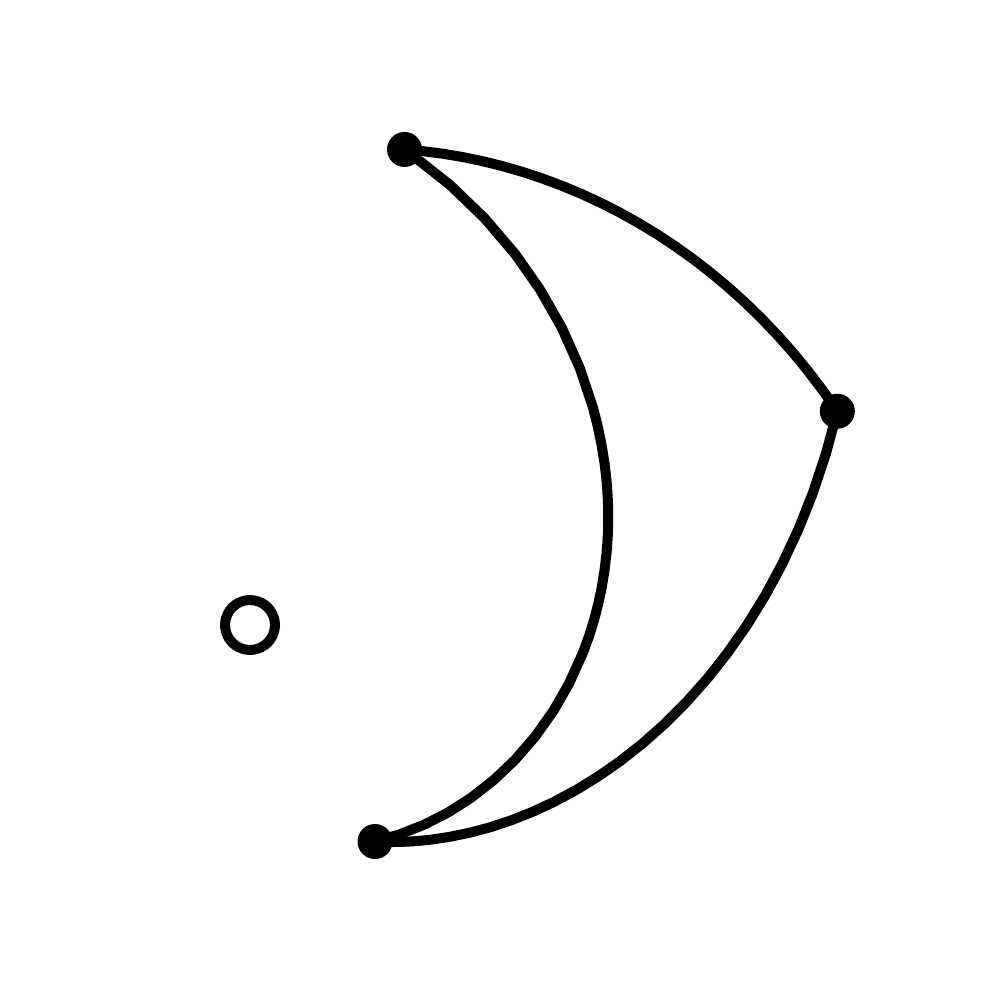}
    \includegraphics[height=5cm]{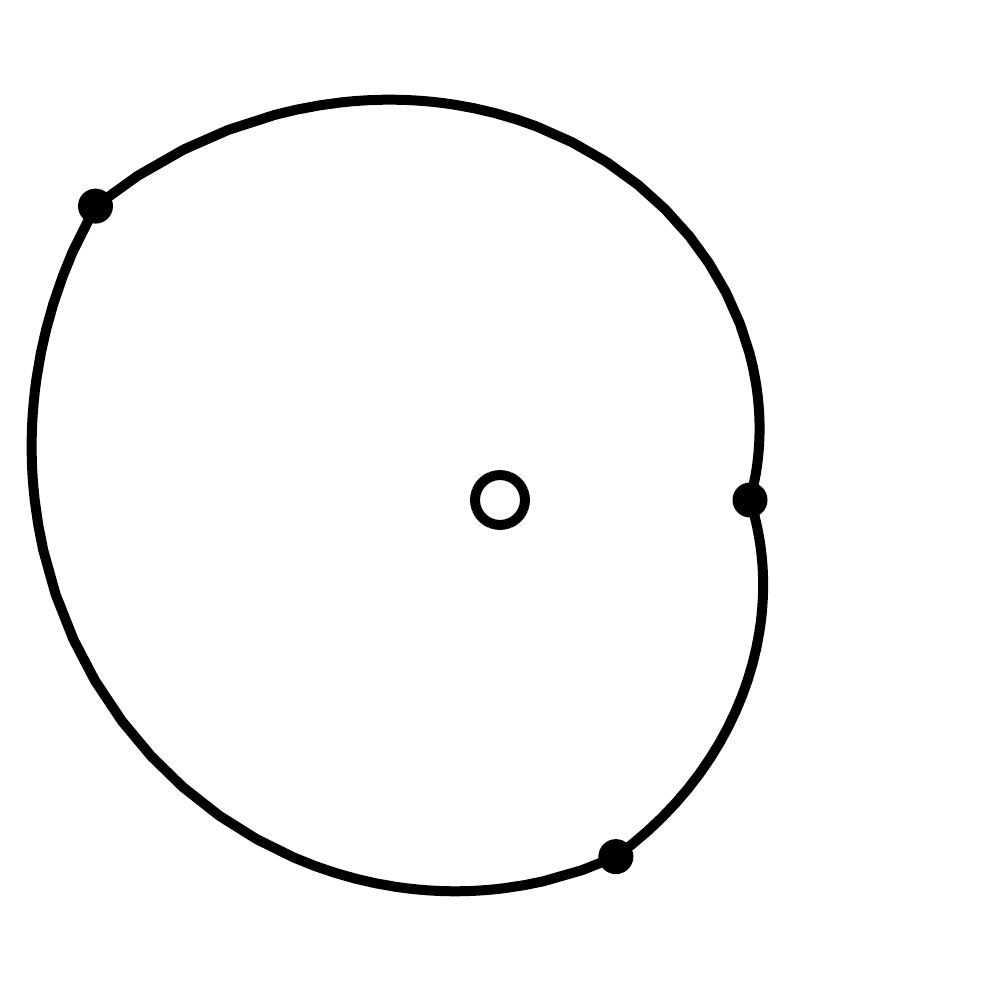}
    \caption{An example of a quasihyperbolic triangle (left) and a quasihyperbolic trigon (right). The small circle stands for the puncture at the origin.
\label{qhtris} \cite{k3}}
  \end{center}
\end{figure}


\subsection{Theorem}\label{heron} (\cite{k3})
  Let $\trik{x}{y}{z}$ be a quasihyperbolic triangle. Then the quasihyperbolic area of $\trik{x}{y}{z}$ is
  $$
    \sqrt{s\big( s-k(x,y) \big) \big( s-k(y,z) \big) \big( s-k(z,x) \big) },
  $$
  where $s = \big( k(x,y)+k(y,z)+k(z,x) \big) /2$.

It is a natural question to ask whether for some other domains
the Law of Cosines holds as an inequality, see \cite{k3}. For the
case of the half plane the problem was solved in \cite{HPWW}.

\renewcommand{\H}{\mathbb{H}}

\subsection{Lemma}\label{inequalityofcosinH} (\cite{HPWW})
  Let $x,y,z \in \H^2$ be distinct points. Then
  $$
    k_{\H^2}(x,y)^2 \ge k_{\H^2}(x,z)^2+k_{\H^2}(y,z)^2-2 k_{\H^2}(y,z)k_{\H^2}(x,z) \cos \gamma,
  $$
  where $\gamma$ is the Euclidean angle between geodesics $J_k[z,x]$ and $J_k[z,y]$.



These results raise many questions about generalizations to more
general situations. For instance, what about the case of domains
with finitely many boundary points or simple domains such as a
sector, a strip or a polygon?


\end{document}